\documentclass[letterpaper]{siamart171218}

\usepackage{amsfonts,amsmath,amssymb}
\usepackage{graphicx}

\newcommand{\bs}[1]{{\boldsymbol{#1}}}

\begin{document}

\title{A note on the relationship between PDE-based precision operators and Mat\'ern covariances}
\author{Umberto Villa and Thomas O'Leary-Roseberry 
\thanks{Umberto Villa and Thomas O'Leary-Roseberry are with the Oden Institute for Computational Engineering and Sciences, University of Texas at Austin, Austin, TX 78712.}
\thanks{Further author information: (Send correspondence to Umberto Villa.)\\ E-mail: uvilla@oden.utexas.edu, Telephone: +1 512 232-3453} }
\maketitle

\begin{abstract}
The purpose of this technical note is to summarize the relationship between the marginal variance and correlation length of a Gaussian random field with Mat\'ern covariance and the coefficients of the corresponding partial-differential-equation (PDE)-based precision operator.
\end{abstract}

\section{Introduction}
The seminal work of Lindgren et al \cite{LindgrenRueLindstrom11} establishes a link between covariance operators of the Mat\'ern family and precision operators defined in terms of partial differential equations (PDEs).

In particular, we consider a Gaussian random field of the Mat\'ern family defined over $\Omega \equiv \mathbb{R}^d$ with Mat\'ern covariance function defined by

\begin{equation}
\label{eq:matern_cov}
c(\bs{r}) = \frac{\sigma^2}{2^{\nu -1} \Gamma(\nu)}\left( \kappa \| \bs{r}\|\right)^\nu K_\nu\left( \kappa \|\bs{r}\|\right).
\end{equation}

Above, $\bs{r} = \bs{x}_1 - \bs{x}_2$ is the distance between any two points $\bs{x}_1, \bs{x}_2 \in \mathbb{R}^d$, $\sigma^2$ is the marginal variance, $\kappa > 0$ is a scale parameter (inversely proportional to the correlation length), $\Gamma$ is the gamma function, and $K_\nu$ is the modified Bessel function of the the second kind and order $\nu > 0$.  The parameter $\nu$ is the smoothness parameter of the field. Empirically, we define the correlation length $\rho$ as
\begin{equation}
\label{eq:corr_len}
    \rho = \frac{\sqrt{8\nu}}{\kappa},
\end{equation}
which corresponds to correlations near 0.1 at the distance $\rho$.

As shown in \cite{LindgrenRueLindstrom11}, a Gaussian random field $u(\bs{x})$ with zero mean and Mat\'ern covariance given by \eqref{eq:matern_cov} is the solution of the linear fractional stochastic PDE

\begin{equation}
\label{eq:spde}
    \left( \kappa^2 - \Delta\right)^{\alpha/2} u = s \mathcal{W}, \quad \forall \bs{x} \in \Omega \equiv \mathbb{R}^d,
\end{equation}
where $\mathcal{W}$ denotes white noise in $\Omega$, and the exponent $\alpha$ is given by
\begin{equation}
 \alpha = \nu + \frac{d}{2}.
 \end{equation}
 
Above, the scaling parameter $s$ is chosen such that $u(\bs{x})$ has marginal variance $\sigma^2$. In particular, we have \cite{LindgrenRueLindstrom11}
 
 \begin{equation}
 \label{eq:scaling_factor}
 s = \sigma\kappa^\nu \sqrt{ \frac{\Gamma(\nu + d/2)(4\pi)^{d/2}}{\Gamma(\nu)}}.
 \end{equation}
 
 \section{The  biLaplacian prior precision in hIPPYlib}
In this section, we explain how to derive the coefficients $\gamma$ and $\delta$ of the PDE-based precision operator $\mathcal{A}^2$ with
\begin{equation}
\label{eq:sqrt_precision}
\mathcal{A} = \delta I - \nabla \cdot (\gamma \nabla ),
\end{equation}
implemented in the hIPPYlib library\cite{VillaPetraGhattas21,VillaPetraGhattas18}.

\subsection{The two-dimensional case}
\label{sec:2D}
 For the case of the biLaplacian prior in 2D, we assume
 \begin{equation}
 \label{eq:2dvals}
     d = 2, \quad \nu = 1,
 \end{equation}
so that the exponent $\alpha/2 = 1$ in \eqref{eq:spde} is an integer value. Substituting \eqref{eq:2dvals} into \eqref{eq:scaling_factor}, the expression of the scaling factor simplifies to
\begin{equation}
\label{eq:scalingBiLapl2D}
    s_{\rm BiLapl2D} = \sigma\kappa \sqrt{ \frac{\Gamma(2)(4\pi)}{\Gamma(1)}} = 2\sigma\kappa \sqrt{\pi}
\end{equation}
where we used the fact that $\Gamma(2) = (2-1)! = 1$ and $\Gamma(1) = (1-1)! = 1$.

Thus we can rewrite \eqref{eq:spde} as

\begin{equation}
    \label{eq:spde_BiLapl2D}
    \left( \kappa^2 - \Delta\right) u = 2\sigma\kappa \sqrt{\pi} \mathcal{W}\quad \forall \bs{x} \in \Omega \equiv \mathbb{R}^2.
\end{equation}

Finally, by dividing both sides of \eqref{eq:spde_BiLapl2D} by the scaling factor $s$ in \eqref{eq:scalingBiLapl2D} we rewrite \eqref{eq:spde_BiLapl2D} in the form used by hIPPYlib as
\begin{equation}
\label{eq:hippylib_biLapl2D}
    \left( \delta - \gamma \Delta\right) u = \mathcal{W} \quad \forall \bs{x} \in \Omega \equiv \mathbb{R}^2,
\end{equation}
where the coefficients $\gamma$ and $\delta$ are defined as
\begin{equation}
\label{eq:gamma2D}
    \gamma = \frac{1}{2\sigma\kappa \sqrt{\pi}} = \frac{1}{4\sqrt{2\pi}}\frac{\rho}{\sigma},
\end{equation}
and
\begin{equation}
\label{eq:delta2D}
    \delta = \frac{\kappa^2}{2\sigma\kappa \sqrt{\pi}} = \frac{1}{2\sqrt{\pi}}\frac{\kappa}{\sigma} = \frac{\sqrt{2}}{\sqrt{\pi}}\frac{1}{\sigma\rho},
\end{equation}
To obtain the expression above with used the fact that, when $\nu=1$, \eqref{eq:corr_len} gives $\kappa = \sqrt{8}/\rho$.

\subsection{The three-dimensional case}
For the case of the biLaplacian prior in 3D, we assume 
\begin{equation}
 \label{eq:3dvals}
     d = 3, \quad \nu = \frac{1}{2},
 \end{equation}
so that the exponent $\alpha/2 = 1$ in \eqref{eq:spde} is an integer value. Note that this choice of the smoothness parameter $\nu$ corresponds to the special case in which the Mat\'ern covariance function \eqref{eq:matern_cov} reduces to the exponential covariance function
\begin{equation*}
c(\bs{r}) = \sigma^2 \exp\left( - \kappa \| \bs{r} \| \right).
\end{equation*}

Substituting \eqref{eq:3dvals} into \eqref{eq:scaling_factor} and \eqref{eq:corr_len}, recalling that $\Gamma(1/2) = \sqrt{2}$, we obtain
\begin{equation}
\label{eq:scalingBiLapl3D}
    s_{\rm BiLapl3D} = \sigma \sqrt{\kappa} \sqrt{ \frac{\Gamma(2)(4\pi)^{3/2}}{\Gamma(\frac{1}{2})}} = 2\sigma \sqrt{\kappa} \sqrt{2\pi},
\end{equation}
and
\begin{equation}
    \kappa = \frac{2}{\rho}.
\end{equation}

Following the same calculation as is Section \ref{sec:2D}, we obtain the following expression for $\gamma$ and $\delta$:
\begin{equation}
\label{eq:gamma3D}
    \gamma = \frac{1}{2\sigma \sqrt{\kappa} \sqrt{2\pi}} = \frac{1}{4\sqrt{\pi}}\frac{\sqrt{\rho}}{\sigma},
\end{equation}
and
\begin{equation}
\label{eq:delta3D}
    \delta = \frac{\kappa^2}{2\sigma \sqrt{\kappa} \sqrt{2\pi}} = \frac{1}{2\sqrt{2\pi}} \frac{\kappa^{3/2}}{\sigma} = \frac{1}{\sqrt{\pi}}\frac{1}{\sigma \rho^{3/2}}.
\end{equation}

\subsection{Extension to spatially varying coefficients and bounded domains}
Note that the equivalence between \eqref{eq:matern_cov} and the precision operator $\mathcal{A}^2$ with $\mathcal{A}$ given in \eqref{eq:sqrt_precision} only holds for infinite domains and constant coefficients $\gamma$ and $\delta$.

When the domain $\Omega$ is bounded with well-defined unit normal vector $\bs{n}$ to the boundary $\partial\Omega$, the PDE-based precision operator $\mathcal{A}^2$ is defined in hIPPYlib with

\begin{equation}
    \label{eq:isotropic_pde}
    \mathcal{A} u = 
    \left\{
    \begin{array}{ll}
    \left( \delta - \nabla \cdot (\gamma \nabla ) \right) u &  \forall \bs{x} \in \Omega \subset \mathbb{R}^d\\
    \gamma \nabla u \cdot \bs{n} + \beta u  & \forall \bs{x} \in \partial\Omega,
    \end{array}
    \right.
\end{equation}

where the coefficient $\beta = \frac{\sqrt{\delta\gamma}}{1.42}$ is empirically chosen as in \cite{Daon} to reduce boundary artifacts, i.e., the inflation of the marginal variance near $\partial\Omega$.

Additionally, one can introduce spatial anisotropy in the samples with the introduction of the symmetric positive-definite tensor $\Theta \in \mathbb{R}^{d\times d }$ in the definition of the PDE-based precision operator:
\begin{equation}
    \label{eq:anisotropic_pde}
    \mathcal{A}_\Theta u = 
    \left\{
    \begin{array}{ll}
    \left( \delta - \nabla \cdot (\gamma \Theta \nabla ) \right) u &  \forall \bs{x} \in \Omega \subset \mathbb{R}^d\\
    \gamma \nabla u \cdot \bs{n} + \beta u  & \forall \bs{x} \in \partial\Omega,
    \end{array}
    \right.
\end{equation}
The anisotropic tensor $\Theta$ can be any user-prescribed symmetric positive definite matrix, and will introduce directional structure based on its eigenvectors. In hIPPYlib, rotation matrices are implemented. For example in 2D, hIPPYlib implements $\Theta$ of the form
\begin{equation*}
    \Theta = \begin{bmatrix}
        \theta_1\sin(\alpha)^2 & (\theta_1-\theta_2)\sin(\alpha)\cos(\alpha)\\
        (\theta_1-\theta_2)\sin(\alpha)\cos(\alpha) & \theta_2\cos(\alpha)^2
    \end{bmatrix},
\end{equation*}
where $\theta_1,\theta_2,\alpha$ are user-prescribed parameters that control the spatial orientation of the resulting anisotropic correlation structure.

\section{Numerical Examples in 1D and 2D}

In this section we demonstrate some numerical examples of the BiLaplacian prior implemented in hIPPYlib \cite{VillaPetraGhattas21,VillaPetraGhattas18}, see notebook \texttt{7\_GaussianPriors.ipynb} therein for reference. We focus on 1D and 2D unit domains $\Omega = (0,1)^d$. For all examples we utilize piecewise linear finite elements. 

\subsection{Isotropic Gaussian Random Fields}

\begin{figure}[h]
    \centering
    \includegraphics[width=0.9\textwidth]{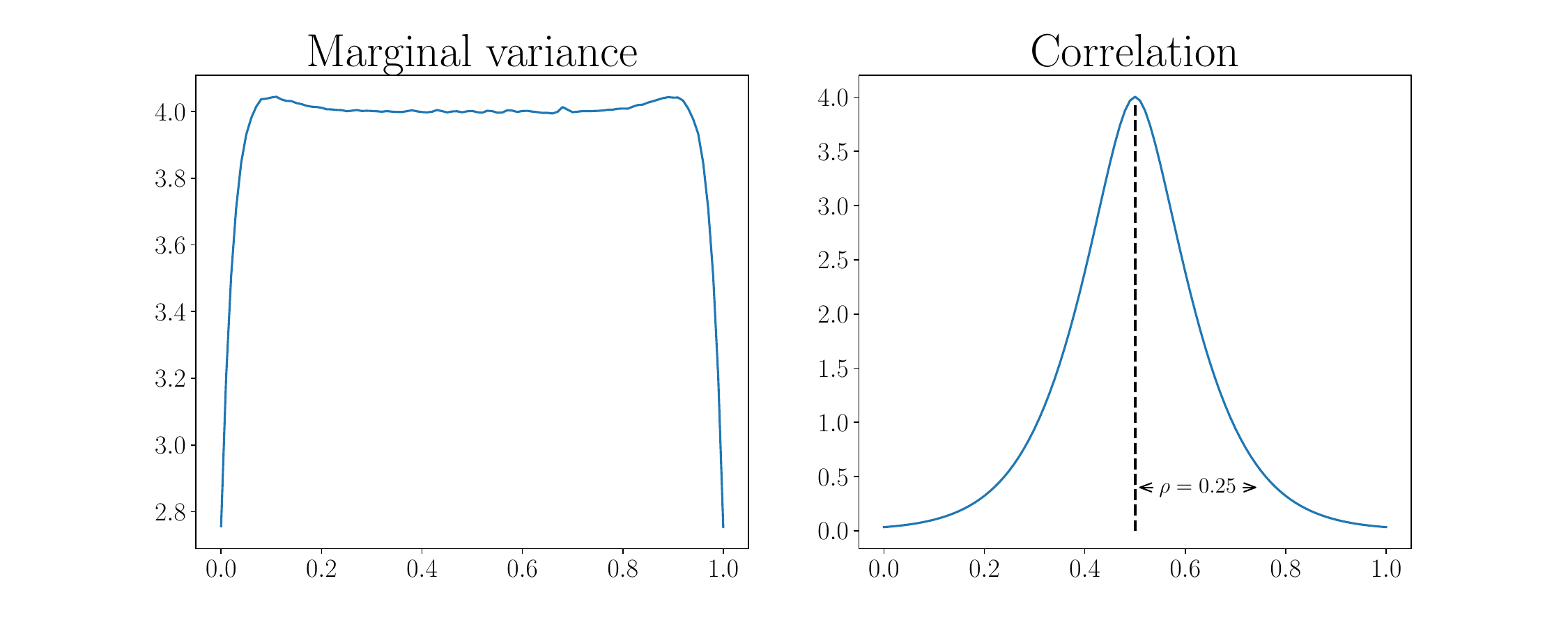}
    \caption{Marginal variance (left) and correlation length structure (right) for the biLaplacian prior with Robin boundary conditions in 1D. On the left one can see the effects of the boundary on the marginal variance; in the interior of the domain the marginal variance is close to $4$, but near the boundary it decreases to $\sim2.8$. On the right, the correlation is shown at the point $x =0.5$; the effects of the boundary are limited as the correlation decreases by a factor of $10$ a distance of $\rho =0.25$ away.}
    \label{fig:1d_robin_var_corr}
\end{figure}

We begin with studying 1D isotropic Gaussian random fields as in \eqref{eq:isotropic_pde} to investigate the correlation and marginal variance and how they are effected by the use of boundary conditions. For this set of examples we set the marginal variance $\sigma^2 = 4$, and correlation length $\rho = 0.25$. We consider the unit segment domain discretized using a mesh with $n_x=100$ elements. First in Figure \ref{fig:1d_robin_var_corr} we demonstrate the marginal variance and correlation structure when utilizing the Robin boundary condition in \eqref{eq:isotropic_pde}. In this figure the effects of the boundary are limited; the marginal variance is very close to $\sigma^2=4$ throughout the domain $\Omega = (0,1)$, and the correlation length roughly corresponds to a $10\times$ reduction in correlation at a distance $\rho$ away from the correlating collocation point. 

\begin{figure}[h]
    \centering
    \includegraphics[width=0.9\textwidth]{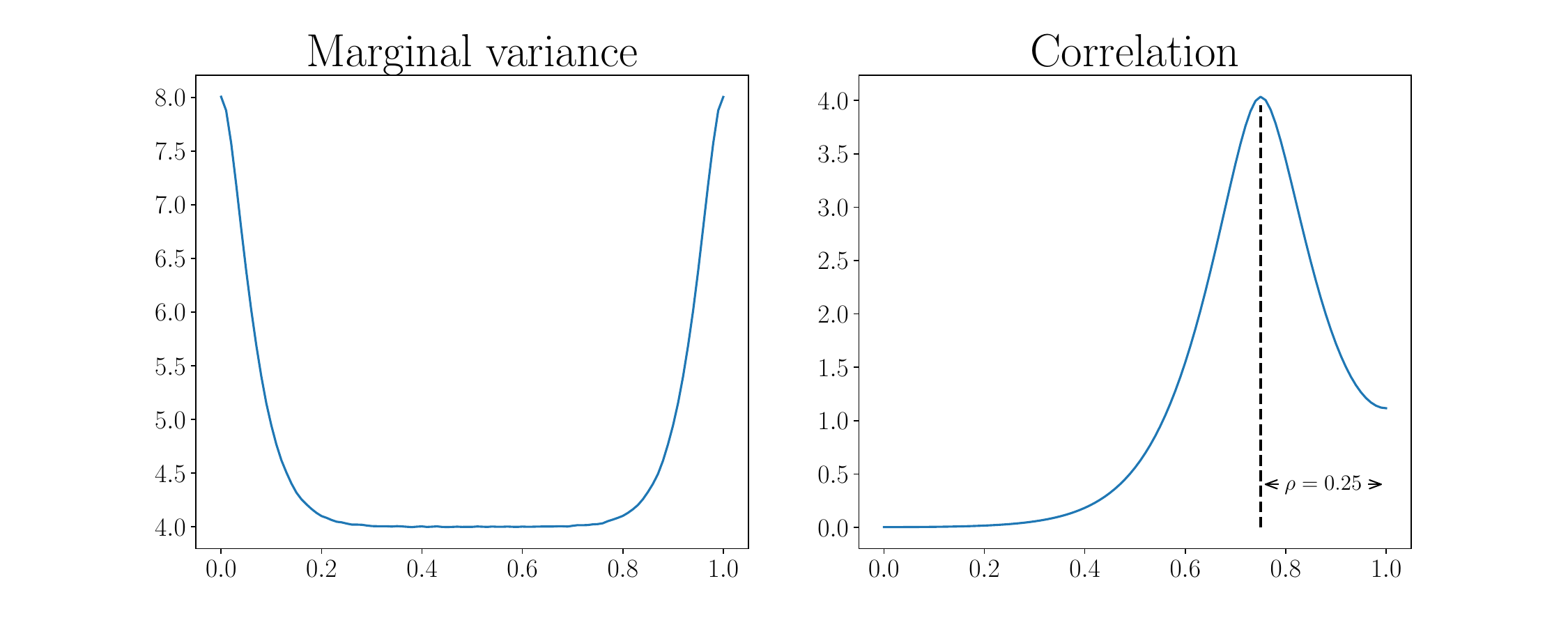}
    \caption{Marginal variance (left) and correlation length structure (right) for the biLaplacian prior with homogeneous Neumann boundary conditions in 1D. On the left one can see the effects of the boundary on the marginal variance; in the interior of the domain the marginal variance is close to $4$, but near the boundary it increases to $\sim8$. Additionally the marginal variance in the interior is more substantially polluted in the interior than in Figure \ref{fig:1d_robin_var_corr}. On the right, the correlation is shown at the point $x =0.75=1.0 - \rho$. In this instance the boundary effects substantially effect the correlation structure near to the boundary.}
    \label{fig:1d_neumann_var_corr}
\end{figure}

In Figure \ref{fig:1d_neumann_var_corr}, we demonstrate the boundary artefacts that are amplified if one utilizes an homogeneous Neumann boundary condition in \eqref{eq:isotropic_pde}, instead of the Robin BC. Indeed, in this case the marginal variance increases substantially at the boundary. Additionally this figure demonstrates how this pollutes the correlation length at a distance $\rho=0.25$ away from the boundary, as in this case the correlation no longer is reduced by a factor of $10$, instead it is closer to a reduction of $2.5$.

\begin{figure}[h]
    \centering
    \includegraphics[width=0.9\textwidth]{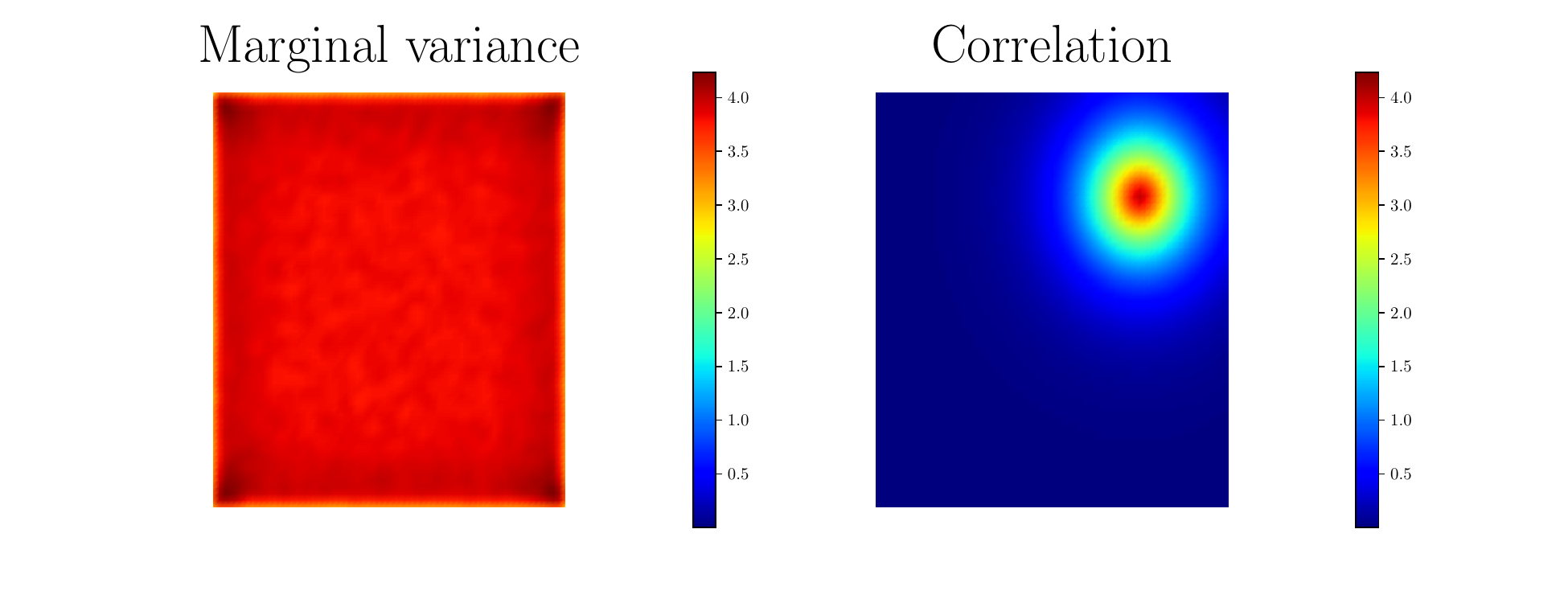}
    \caption{Marginal variance (left) and correlation length structure (right) for the biLaplacian prior with Robin boundary conditions in 2D. On the left one can see the effects of the boundary on the marginal variance; in the interior of the domain the marginal variance is close to $4$, but near the boundary it is reduced a small amount. On the right, the correlation is shown at the point $\mathbf{x} = (0.75,0.75) = (1-\rho,1-\rho)$; the effects of the boundary are limited as the correlation decreases by a factor of $10$ a distance of $\rho =0.25$ away.}
    \label{fig:2d_robin_var_corr}
\end{figure}

\begin{figure}[h]
    \centering
    \includegraphics[width=0.9\textwidth]{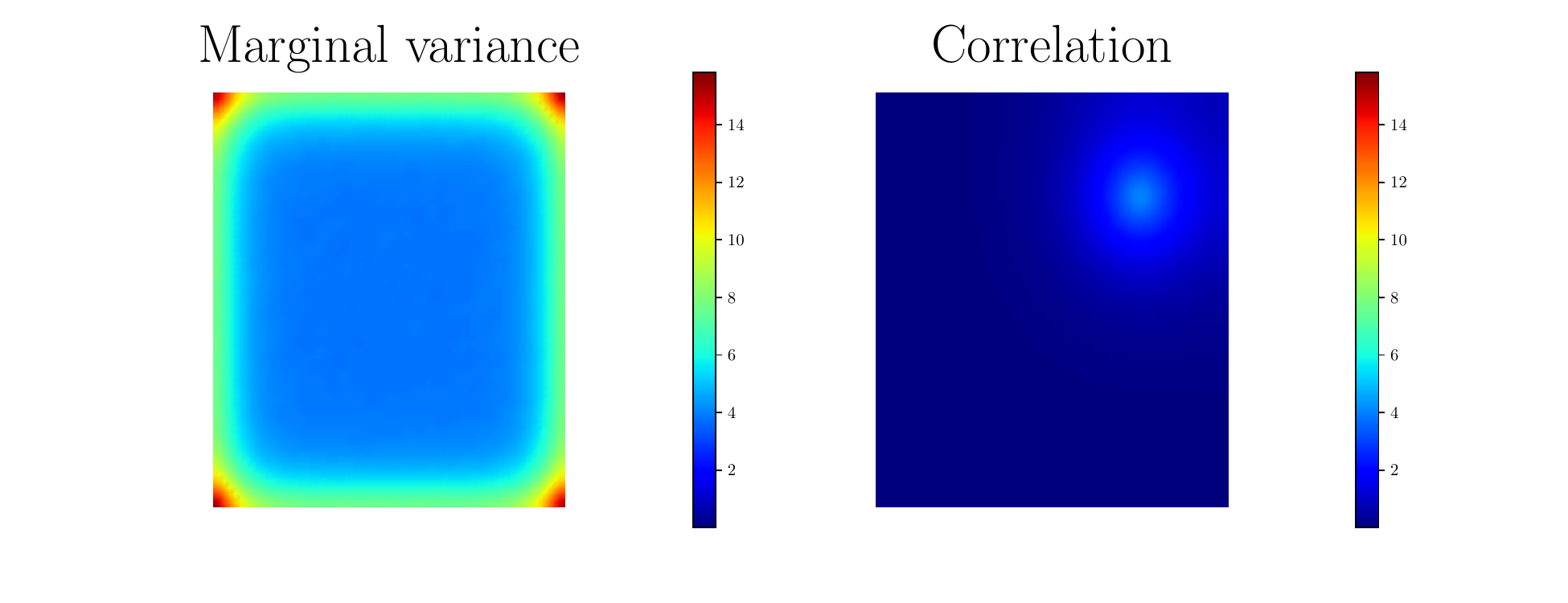}
    \caption{Marginal variance (left) and correlation length structure (right) for the biLaplacian prior with homogeneous Neumann boundary conditions in 2D. On the left one can see the effects of the boundary on the marginal variance; in the interior of the domain the marginal variance is close to $4$, but near the boundary it increases, reaching $\sim15$ at the corners. Additionally the marginal variance in the interior is more substantially polluted in the interior than in Figure \ref{fig:2d_robin_var_corr}. On the right, the correlation is shown at the point $\mathbf{x} = (0.75,0.75) = (1-\rho,1-\rho)$. In this instance the boundary effects substantially effect the correlation structure near to the boundary.}
    \label{fig:2d_neumann_var_corr}
\end{figure}

We proceed with investigation of the marginal variance and correlation structure for the same isotropic biLaplacian prior in 2D, and with $\sigma^2 = 4$ and $\rho = 0.25$. For these examples we consider a unit square domain and use a mesh with $(n_x,n_y) = (64,64)$ elements. In Figure \ref{fig:2d_robin_var_corr} we show the marginal variance and correlation structure for the isotropic prior corresponding to \eqref{eq:isotropic_pde} with Robin boundary conditions. In Figure \ref{fig:2d_neumann_var_corr} we show the marginal variance and correlation structure for the same prior but with homogeneous Neumann boundary conditions. In both cases we demonstrate the correlation structure for the point $\mathbf{x} = (0.75,0.75) = (1-\rho,1-\rho)$, demonstrating the substantial boundary artifacts that arise when the homogeneous Neumann boundary condition is used. 

\begin{figure}[h]
    \centering
    \includegraphics[width=0.9\textwidth]{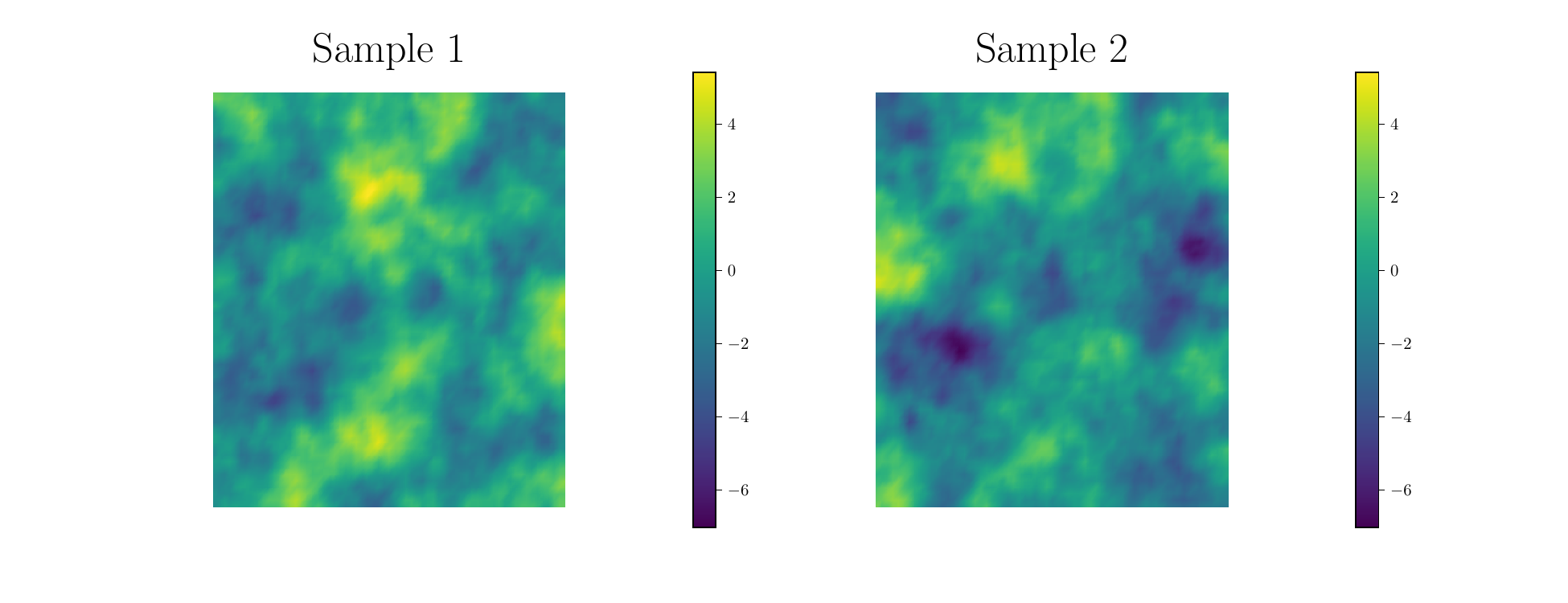}
    \caption{Two isotropic samples from the biLaplcian prior with Robin boundary conditions on a $64\times 64$ mesh.}
    \label{fig:2d_isotropic samples}
\end{figure}

In Figure \ref{fig:2d_isotropic samples} we show two instances of samples from the biLaplacian prior \eqref{eq:isotropic_pde} utilizing Robin boundary conditions.

\newpage
\subsection{Anisotropic Gaussian Random Fields}

\begin{figure}[h]
    \centering
    \includegraphics[width=1.0\textwidth]{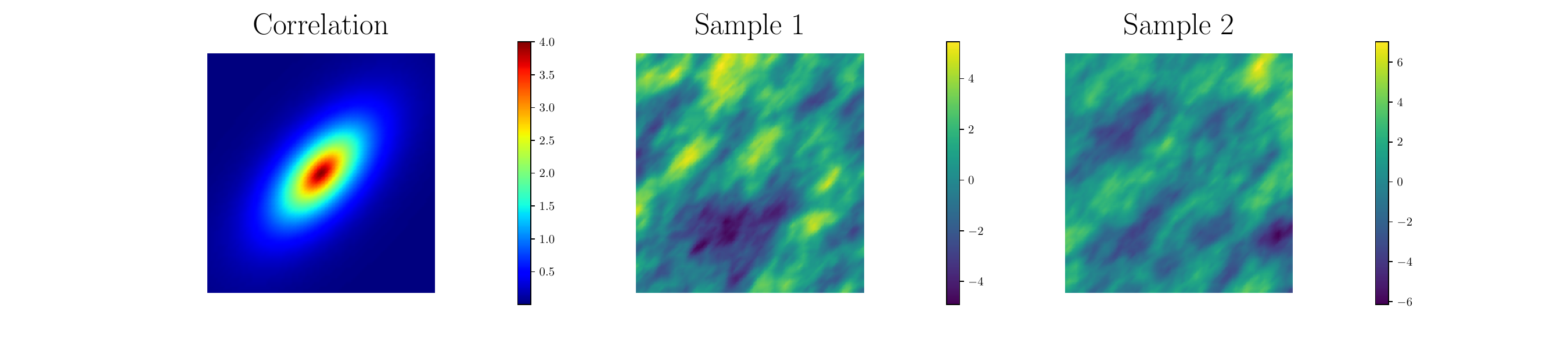}
    \includegraphics[width=1.0\textwidth]{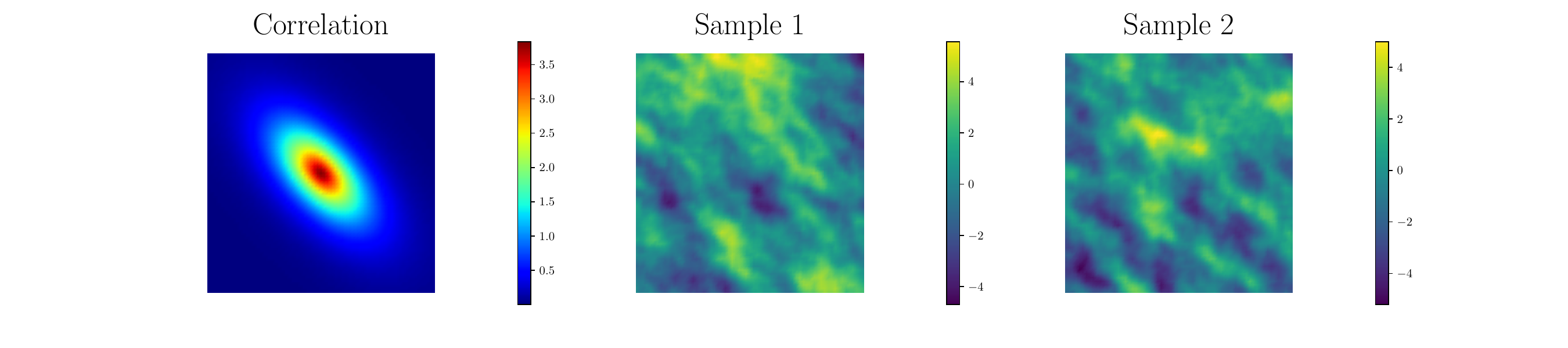}
    \caption{The top row shows the correlation structure and two samples from the distribution with $\alpha= \frac{\pi}{4}$. The bottom row shows the correlation structure and two samples from the distribution with $\alpha= -\frac{\pi}{4}$. The correlation structure is rotated by angle $\alpha$}
    \label{fig:anisotropic_corr_var_pi_4}
\end{figure}

Next we investigate the effects of the anisotropic tensor $\Theta$ utilized in \eqref{eq:anisotropic_pde}. We use the same numerical setup as before, and utilize the Robin boundary conditions. We set $\Theta$ are $\theta_1 = 2,\theta_2 = 0.5$, and we vary $\alpha = \pm\frac{\pi}{4}$ to demonstrate the effects of rotating the orientation of the correlation structure; see Figure \ref{fig:anisotropic_corr_var_pi_4}.

\newpage

\bibliographystyle{siam}

\bibliography{references,local}

\end{document}